\documentclass{siamltex}
\usepackage{xcolor} 
\usepackage{amsfonts}
\usepackage{amssymb}
\usepackage{amsfonts}
\usepackage{mathtools}
\usepackage{enumitem}
\usepackage{color}

\newcommand{\R}{{\mathbb R}}
\newcommand{\cE}{{\mathcal{E}}}
\newcommand{\cG}{{\mathcal{G}}}
\newcommand{\cW}{{\mathcal{W}}}
\newcommand{\cN}{{\mathcal{N}}}

\newtheorem{rmk}{Remark}

\newtheorem{exmp}{Example}[section]
\newcommand{\bl}{\begin{list}{ \ }{
\leftmargin=.325in}}

\title{Estimating and increasing the structural robustness of a network}\medskip
\author{
Silvia Noschese\thanks{Dipartimento di Matematica ``Guido Castelnuovo'', SAPIENZA 
Universit\`a di Roma, P.le A. Moro, 2, I-00185 Roma, Italy. E-mail: 
{\tt noschese@mat.uniroma1.it}}
\and
Lothar Reichel\thanks{Department of Mathematical Sciences, Kent State University, Kent,
OH 44242, USA. E-mail: {\tt reichel@math.kent.edu}}
}
\begin{document}

\maketitle

\begin{abstract}
The capability of a network to cope with threats and survive attacks is referred to as
its robustness. This paper discusses one kind of robustness, commonly denoted structural 
robustness, which increases when the spectral radius of the adjacency matrix associated 
with the network decreases. We discuss computational techniques for identifying edges,
whose removal may significantly reduce the spectral radius. Nonsymmetric adjacency 
matrices are studied with the aid of their pseudospectra. In particular, we consider
nonsymmetric adjacency matrices that arise when people seek to avoid being infected by 
Covid-19 by wearing facial masks of different qualities.
\end{abstract}


\section{Introduction}\label{sec1}
Networks appear in many areas, including transportation, communication, social science, 
and chemistry; see, e.g., Estrada \cite{Es} and Newman \cite{Ne} for many examples. An 
edge-weighted network is represented by a graph $\cG=\{\cN,\cE,\cW\}$, which consists of a 
set of \emph{nodes} $\cN=\{n_j\}_{j=1}^n$, a set of \emph{edges} $\cE=\{e_j\}_{j=1}^m$ 
that connect the nodes, and a set of \emph{edge weights} $\cW=\{w_j\}_{j=1}^m$ that 
indicate the importance of the edges. The weights are assumed to be positive. For 
instance, in a road network, the nodes $n_j$ may represent cities, the edges $e_j$ may 
represent roads between the cities, and the edge weight $w_j$ may be proportional to the 
amount of traffic on the road represented by edge $e_j$. We refer to a graph $\cG$ as 
\emph{undirected} if for each edge $e_j$, there is an edge $e_{k_j}$ that points in the 
opposite direction and has the same weight as $e_j$. If this is not the case, then the 
graph $\cG$ is said to be \emph{directed}. 

The \emph{adjacency matrix} $A=[a_{ij}]_{i,j=1}^n\in{\R}^{n\times n}$ associated with the 
graph $\cG$ has the entry $a_{ij}=w_k$ if there is an edge $e_k$ emerging from node $n_i$ 
and ending at node $n_j$; if the graph is undirected, then also $a_{ji}=w_k$. Other matrix
entries vanish. Thus, the matrix $A$ is symmetric if and only if the graph $\cG$ is 
undirected. We will assume that there are no self-loops and no multiple edges. The former
implies that the diagonal entries of $A$ vanish. Typically, the number of edges, $m$,
satisfies $1\le m\ll n^2$. Then the matrix $A$ is sparse. 

The maximum of the magnitudes of the eigenvalues of $A$ is known as the 
\emph{spectral radius} of $A$. We will denote the spectral radius of $A$ by $\rho(A)$. It 
has been shown that the spectral radius is an important indicator of how flu-type 
infections spread in the network that is associated with the adjacency matrix $A$; the
smaller $\rho(A)$, the less spread; see, e.g., \cite{JKVV,MSN} and below. This paper seeks
to shed light on how the spectral radius of an adjacency matrix can be reduced by targeted 
edge perturbations, i.e., by reducing edge-weights or removing edges. It is well known 
that reducing an edge-weight, or removing an edge, does not increase the spectral radius 
of a nonnegative matrix; see, e.g., \cite[Corollary 8.1.19]{HJ}. We are interested in 
identifying which weights should be reduced, or which edges should be removed, to achieve
a significant decrease of the spectral radius.

Howard et al. \cite{Ho} discuss the benefits of wearing facial masks to reduce Covid-19
transmission. Several studies found $70\%$ or higher efficacy of facial masks in 
protecting the wearer of Covid-19 infections. They found that wearing a mask protects 
people around persons wearing masks, as well as the people who wear a mask, but to lesser 
degree. Also the type of mask is important; see also Gandhi et al. \cite{GM,GR} for 
related discussions.  

Let the nodes in a graph represent people and the edge weights represent the possibility 
of getting a sufficient viral load to become ill with Covid-19. The modeling of facial 
masks of different quality results in a nonsymmetric adjacency matrix $A$ associated with 
the graph; see Section \ref{sec3}. We refer to a network (or graph) that is robust against
the spread of viruses as \emph{structurally robust}, and measure the structural robustness
in terms of the spectral radius of the adjacency matrix associated with the graph. A 
network is more structurally robust, the smaller the spectral radius is. We are interested 
in determining which edge-weights should be reduced, or which edges should be removed, to
give a relatively large decrease in the spectral radius of the adjacency matrix.

This paper is organized as follows. Section \ref{sec2} discusses the structural robustness
of a network. In particular, the sensitivity of the eigenvalues to perturbations of the
adjacency matrix $A$ associated with the network is considered. Tridiagonal adjacency matrices that 
model the role of face masks are described. Section \ref{sec3} is concerned with the 
calculation of the spectral radius of a large matrix, and with the determination of edges 
that should be eliminated, or whose weight should be reduced, to achieve a relatively
large decrease in the spectral radius of $A$. Properties of the pseudospectrum of a matrix
are reviewed and applied. Some large-scale computed examples are presented in Section 
\ref{sec4}, and concluding remarks can be found in Section \ref{sec5}.

\section{Structural robustness}\label{sec2}
A formulation of structural robustness against the spread of viruses comes from 
spectral graph theory. Epidemiological 
theory predicts that if the effective infection rate of a virus in an epidemic is below 
the reciprocal of the spectral radius $\rho(A)$ of the adjacency matrix $A$ associated 
with the graph that represents the network, then the virus contamination in the network 
dies out over time. In more detail, assume a universal virus birth rate $\beta$ along 
each edge that is connected to an infected node, and a virus death rate $\delta$ for 
each infected node. If the effective infection rate, given by $\beta/\delta$, is 
below the epidemic threshold for the network, i.e., if
\[
\frac{\beta}{\delta}<\frac{1}{\rho(A)},
\]
then the infections tend to zero exponentially over time. In fact, the reciprocal of
the spectral radius $\rho(A)$ is a network-epidemic threshold in a 
Susceptible-Infectious-Susceptible (SIS) network, in which the evolution of the viral 
state $s_i(t)\geq 0$ of node $n_i$, $i=1,2,\ldots,n$, at time $t$ is governed by the 
system of differential equations
$$
\frac{d\mathbf{s}(t)}{dt}=-\delta \mathbf{s}(t) + \beta 
\diag(\mathbf{e}-\mathbf{s}(t))A\mathbf{s}(t),
$$
where $\mathbf{s}(t)=[s_1(t),s_2(t),\ldots,s_n(t)]^T$ and 
$\mathbf{e}=[1,1,\ldots,1]^T\in\R^n$. Indeed, if one has $\rho(\delta^{-1}\beta A)<1$, then 
$\mathbf{s}(t)\rightarrow\mathbf{0}$ as $t \rightarrow \infty$; see, e.g., \cite{JKVV}
and references therein. The smaller $\rho(A)$ is, the higher is the structural robustness 
of the network against the spread a virus. Hence, in order to enhance the structural 
robustness of a network, one may want to reduce the weights of suitable edges in $\cE$ of
the graph $\cG$, or eliminate certain edges; see \cite{MSN}.

Let $\mathbf{e}_j=[0,\ldots,0,1,0,\ldots,0]^T\in\R^n$ denote the $j$th axis vector and 
assume that the entry $a_{hk}$ of $A$ is positive. Consider the rank-1 matrix 
\[
E_{hk}=-a_{hk}\mathbf{e}_h\mathbf{e}_k^T,
\]
where the superscript $^T$ denotes transposition. Regard the perturbed adjacency matrix
\begin{equation}\label{Atilde}
\widetilde{A}=A+\varepsilon E_{hk},
\end{equation}
where $\varepsilon>0$ is chosen small enough so that the matrix $\widetilde{A}$ is
nonnegative. Assume for the moment that the graph $\cG$ associated with the adjacency 
matrix $A$ is strongly connected, i.e., that starting at any node of the graph, one can 
reach any other node of the graph by traversing the edges along their directions. 
This is equivalent to $A$ being irreducible. Then the Perron--Frobenius theorem 
applies, see, e.g., \cite[Chapter 8]{HJ}, and shows that the eigenvalue of $A$ of largest 
magnitude is unique and equals $\rho(A)$. This eigenvalue is commonly referred to as the 
\emph{Perron root} of $A$. Moreover, right and left eigenvectors of $A$ associated with 
the Perron root,
\[
\mathbf{u}=[u_1,u_2,\ldots,u_n]^T\quad\mbox{and}\quad \mathbf{v}=[v_1,v_2,\ldots,v_n]^T,
\]
respectively, are unique up to scaling. They can be normalized to be of unit Euclidean 
norm and only have positive entries. These normalized vectors are known as the right and 
left \emph{Perron vectors}, respectively, of $A$. We define the \emph{spectral impact} of
the directed edge $e_{hk}\in\cE$ on the spectral radius of $A$ as the relative change of 
the spectral radius $\rho(A)$ induced by the perturbation of the edge \eqref{Atilde}, 
i.e.,
\[
s_{hk}^{(\rho(A))}=\frac{\rho(A)-\rho(\widetilde{A})}{\rho(A)}.
\]

A first order approximation of $s_{hk}^{(\rho(A))}$, when $0<\varepsilon\ll 1$, is derived
in \cite[Eq. (17)]{MSN} as follows. Observe that
\[
\rho(A)-\rho(\widetilde{A})\approx-
\frac{\mathbf{v}^T \varepsilon E_{hk}\mathbf{u}}{\mathbf{v}^T\mathbf{u}}
= \frac{\varepsilon\, a_{hk} v_h u_k}{\mathbf{v}^T\mathbf{u}}>0,
\]
Let $\lambda$ be a simple eigenvalue of $A$, and let $\mathbf{x},\mathbf{y}\in\R^n$
be associated right and left eigenvectors. Then the condition number of $\lambda$ is
defined as 
\[
\kappa(\lambda)=\frac{\|\mathbf{x}\|\|\mathbf{y}\|}{|\mathbf{x}^T\mathbf{y}|};
\]
see, e.g., \cite{NP06}.
Here and throughout this paper $\|\cdot\|$ denotes the Euclidean vector norm or the 
spectral matrix norm. In particular, the condition number of the largest eigenvalue 
$\rho(A)$ of $A$ is given by
\[
\kappa(\rho(A))=\frac{1}{\mathbf{v}^T\mathbf{u}}.
\]
Therefore,
\begin{equation}\label{sp_rad}
s_{hk}^{(\rho(A))}\approx \alpha_{hk} \frac{\varepsilon \, \kappa(\rho(A))}{\rho(A)},
\end{equation}
where 
\begin{equation}\label{alphacoeff}
\alpha_{hk}=a_{hk} v_h u_k.
\end{equation}
Notice that the first order approximation \eqref{sp_rad} of the spectral impact of the 
edge $e_{hk}\in\cE$, which points from node $n_h$ to node $n_k$, depends on the right and left Perron vectors of $A$, as well as on the
weight of the edge $e_{hk}$. To make $\rho(A)$ smaller, we may consider reducing weight(s)
associated with the largest coefficients \eqref{alphacoeff}. To determine these 
coefficients, one needs the Perron vectors $\mathbf{u}$ and $\mathbf{v}$. 
%LRnew 
%Approximations 
%of these vectors for large nonsymmetric adjacency can be conveniently computed by the 
%two-sided Arnoldi method \cite{ZH}; see Section \ref{sec3} for details. 

When the matrix $A$ is symmetric, it is meaningful to require the perturbation of $A$ 
also be symmetric. We therefore define the symmetric perturbation matrix
\[
E^{(S)}_{hk}=-a_{hk}(\mathbf{e}_h\mathbf{e}_k^T+\mathbf{e}_k\mathbf{e}_h^T).
\]
Consider the perturbed matrix 
\[
\widetilde{A}=A+\varepsilon E^{(S)}_{hk}
\]
for some small $\varepsilon>0$. Then a first order approximation of the spectral impact 
$s_{hk}^{(\rho(A))}$ of the undirected edges $e_{hk},e_{kh}\in\cE$ on the spectral radius
$\rho(A)$ is given by
\[
s_{hk}^{(\rho(A))}\approx \alpha_{hk}\frac{\varepsilon}{\rho(A)},
\]
where we have used the fact that the right and left Perron vectors coincide, and
\[
\alpha_{hk}=2\,a_{hk} u_h u_k;
\]
see \cite[Eq. (21)]{MSN}.

\begin{rmk} 
Let the adjacency matrix $A=[a_{hk}]_{h,k=1}^n\in\R^{n\times n}$ be diagonalizable, i.e., 
$A=X\Lambda X^{-1}$, where the columns of $X\in\R^{n\times n}$ are linearly independent 
eigenvectors of $A$, and $\Lambda={\rm diag}[\lambda_1,\lambda_2,\ldots,\lambda_n]$ 
contains the eigenvalues. Then
\begin{equation}\label{rhobds}
\rho(A)^k\leq\|A^k\|\leq \kappa(X)\rho(A)^k,
\end{equation}
where $\kappa(X)=\|X\|\|X^{-1}\|$ is the spectral condition number of $X$. In particular, 
when $A$ is symmetric, we have $\rho(A)^k=\|A^k\|$ for all $k$. 

A \emph{walk} of length $k$ starting at node $n_i$ and ending at node $n_j$ is a sequence 
of $k+1$ nodes $n_{\ell_1},n_{\ell_2},\ldots,n_{\ell_{k+1}}$ with $n_{\ell_1}=n_i$ and 
$n_{\ell_{k+1}}=n_j$ such that there is an edge $e_{q_p}$ that points from node 
$n_{\ell_p}$ to node $n_{\ell_{p+1}}$ for $p=1,2,\ldots,k$; see \cite{Es,Ne}. Edges in a 
walk may be repeated. If the graph is unweighted, then the entry $(i,j)$ of $A^k$ equals 
the number of walks of length $k$ from node $n_i$ to node $n_j$. For weighted graphs, the 
entries of $A^k$ are suitably modified. In view of the bounds \eqref{rhobds}, it may be a 
good idea to eliminate edges in long walks or to reduce the weight of such edges. 

Consider the Frobenius matrix norm $\|A\|_F=\sqrt{\sum_{h,k=1}^n a_{hk}^2}$. The 
inequalities $\|A\|\leq\|A\|_F$ and \eqref{rhobds} suggest that in order to reduce 
$\rho(A)$ the most, it may be a good idea to remove nodes of $\cG$ with many edges, or 
to reduce the weights of edges that emerge from or end at these nodes. In other words, we 
would like to remove or reduce the Euclidean norm of rows and/or columns of the adjacency 
matrix whose Euclidean norm is relatively large. This can we achieved by removing specific
nodes or by reducing edge weights.
\end{rmk}

We conclude this section with a few illustrations for some weighted graphs that are
associated with tridiagonal adjacency matrices. First consider the case when each node 
represents a person, and all persons wear the same kind of facial mask. The persons 
are in a line and each person only can infect the following or preceding person in the 
line. The adjacency matrix is 
\begin{equation}\label{Tmat}
A=\left[ 
\begin{array}{ccccccc}
0 & \sigma &  &  &  &  & \mbox{\Large{O}} \\ 
\sigma & 0& \sigma &  &  &  &  \\ 
& \sigma & 0 & \cdot &  &  &  \\ 
&  & \cdot & \cdot & \cdot &  &  \\ 
&  &  & \cdot & \cdot & \cdot &  \\ 
&  &  &  & \cdot & \cdot & \sigma \\ 
\mbox{\Large{O}} &  &  &  &  & \sigma & 0
\end{array}
\right]\in\R^{n\times n},
\end{equation}
where the edge weight $\sigma>0$ depends on properties of the mask. A high-quality mask 
corresponds to a small value of $\sigma>0$. The graph associated with the matrix 
\eqref{Tmat} is undirected, (strongly) connected, and weighted. 

\begin{proposition}\label{peig}
The Perron root $\rho$ of the nonnegative symmetric tridiagonal Toeplitz matrix 
\eqref{Tmat} is $2\sigma \cos\frac{\pi}{n+1}$. The Perron vector 
$\mathbf{u}=[u_1,u_2,\ldots,u_n]^T$, suitably scaled, has the entries 
$u_k=\sin\frac{k\pi}{n+1}$, $1\leq k\leq n$. In particular, when $n$ is odd, the largest 
entry is $u_{(n+1)/2}$, and when $n$ is even the two largest entries, $u_{n/2}$ and 
$u_{n/2+1}$, have the same size. 
\end{proposition}

\begin{proof}
Explicit formulas for eigenvalues and eigenvectors of tridiagonal Toeplitz matrices can be
found in, e.g., \cite{NPR}. 
\end{proof}

Note that the Perron vector in Proposition \ref{peig} is independent of the numerical
value $\sigma\ne 0$ of the entries of \eqref{Tmat}. Moreover, the Perron vector suggests 
that the node $n_{(n+1)/2}$ for $n$ odd, and the nodes $n_{n/2}$ and $n_{(n+2)/2}$ for $n$
even, are the most important nodes of the graph; see, e.g., Bonacich \cite{Bo}. This is in
agreement with the intuition that the nodes ``in the middle'' of the graph are the best 
connected nodes and, therefore, the most important ones. According to the estimate 
\eqref{sp_rad}, edges that connect these nodes to the graph have the largest spectral 
impact. Consequently, to decrease the spectral radius $\rho(A)$ of the matrix \eqref{Tmat}
maximally, we should reduce the weights of the edges $e_{hk},e_{kh}\in\cE$, where
\begin{itemize}
\item $h=(n+1)/2$ and $k=(n+3)/2$, or $h=(n+1)/2$ and $k=(n-1)/2$, if $n$ is odd;
\item $h=n/2$ and $k=(n+2)/2$ if $n$ is even.
\end{itemize}
Note that setting the edge-weights to zero results in a disconnected graph. It is often 
meaningful to keep a small positive weight. This results in an irreducible adjacency 
matrix. Properties of tridiagonal matrices with some ``tiny'' positive off-diagonal 
entries have been studied by Parlett and V\"omel \cite{PV}. 

\begin{exmp} \label{tridtoep}
Let $A \in \R^{25\times 25}$ be the symmetric tridiagonal Toeplitz matrix \eqref{Tmat} 
with $\sigma=1$. Thanks to Proposition \ref{peig}, one easily computes the spectral 
radius $\rho(A)=1.985418$ and the unit norm Perron vector $\mathbf{u}$. If one chooses to
reduce the weights for the edges $e_{13,14}$ and $e_{14,13}$, as we suggested in the
above discussion, then one obtains the perturbed adjacency matrix for a weighted graph,
\[
\widetilde{A}=A+\varepsilon E^{(S)}_{13,14}= A-\varepsilon(\mathbf{e}_{13}
\mathbf{e}_{14}^T+\mathbf{e}_{14}\mathbf{e}_{13}^T).
\]
Setting $\varepsilon=0.1$ yields $\rho(\widetilde{A})=1.973080$. The spectral impact of 
reducing the weights for the edges $e_{13,14}$ and $e_{14,13}$ is 
$s_{13,14}^{(\rho(A))}=0.006241$, and its first order approximation is 
$2\,\sigma u_{13} u_{14}\varepsilon/\rho(A)=0.007692$.

If, instead, one reduces the weights for the edges $e_{1,2}$ and $e_{2,1}$ and
constructs the perturbed adjacency matrix
\[
\widetilde{A}=A+\varepsilon E^{(S)}_{1,2}= A-\varepsilon(\mathbf{e}_1
\mathbf{e}_2^T+\mathbf{e}_2\mathbf{e}_1^T),
\]
with $\varepsilon=0.1$, one has $\rho(\widetilde{A})=1.985055$. Here
$s_{1,2}^{(\rho(A))}=0.000182$ and $2\,\sigma u_{1} u_{2}\varepsilon/\rho(A) = $ $0.000224$.

Thus, $s_{13,14}^{(\rho(A))}$ can be seen to be significantly larger than 
$s_{1,2}^{(\rho(A))}$. This example shows the reduction of the spectral radius of
the adjacency matrix to be much larger when the weight of an ``important'' edge is 
reduced than when the weight of a less important edge is reduced by the same amount. 
This illustrates the importance of well-connected people wearing high-quality face
masks, which correspond to a small edge weight. We remark that $\varepsilon>0$ is
chosen fairly small in this example so that the estimate \eqref{sp_rad} is applicable.
\end{exmp}

In a description of an epidemic, every node of a realistic network through which an 
infectious disease might spread, may correspond to an individual as well as to a cluster 
of individuals. For instance, we can assign individuals with similar age or location to the 
same group. Such a high-level description may be simpler to analyze 
than a model that accounts for each individual. 
Consider a network in which each node $n_i$ corresponds to a place where a cluster of 
cohabiting people live. This model takes into account the lockdown protocol adopted by the 
Italian Government on Easter 2021, during one of the most delicate phases of the spreading 
of the COVID-19. On Sunday, April 4th, it was permitted to have Easter lunch only with a
small number of friends or relatives.  When meeting people who did not live together, it 
was recommended to wear a face mask indoors, also. In particular, two cohabiting 
adults, possibly accompanied by minor children, were allowed to do only one trip in order 
to reach a place (within the region) of two cohabiting friends or relatives. Moreover, it was 
permitted to host at one place up to two non-cohabiting people, plus minor children.  
As an example, in a family of mother, father, and a 21 years old child living at place $n_i$, 
the parents were permitted to visit two cohabiting relatives or friends at place $n_j$ and the 
child was permitted to receive two cohabiting relatives or friends at place $n_i$. 
Therefore, if a trip corresponds to an edge, then a tridiagonal  matrix can be seen as 
an effective approximation of the 
adjacency matrix of the ``2021 Easter lunch network". Such a simplified model allows one to 
analyze the results and verify that the approach is reasonable at least in this simple 
situation. Moreover, the behavior of such a model may illustrate different scenarios for
decision-makers in public health regarding the lockdown intensity. As an example, taking 
into account aerosol transmission of the COVID-19 virus in enclosed spaces, also the length 
of time spent at a place might be regulated in a hypothetical lockdown protocol to be 
adopted, e.g., for Christmas dinner on December 24th, 2021. This would result in weighting
edges according to the duration of a visit. Alternatively, in order to
mitigate the infectious disease spread, the Italian Government hypothetically might not permit
a family unit or, better, a cluster of cohabiting people to divide.
This would result in a ``2021 Christmas dinner network'' where an edge is weighted according
to the number of people in a moving cluster and a node can have in-edges or 
out-edges.

We now turn to a more accurate model of the role of facial masks. Let node $v_i$ represent 
a person who wears a mask, and assume that the fraction $w_i^{(out)}$ of viruses penetrates
the mask from the outside in unit time, and the fraction $w_i^{(in)}$ penetrates the mask 
from the inside in unit time. Let again the adjacency matrix $A\in\R^{n\times n}$ be 
tridiagonal. The edge from $v_i$ to $v_{i+1}$ has the weight $w_i^{(in)}w_{i+1}^{(out)}$
for $i=1,2,\ldots,n-1$, and the edge from $v_{i+1}$ to $v_i$ has the weight 
$w_{i+1}^{(in)}w_{i}^{(out)}$ for $i=1,2,\ldots,n-1$. This yields the adjacency matrix
{\small
\begin{equation}\label{Amat}
A=\left[ 
\begin{array}{ccccccc}
0 & w_1^{(in)}w_2^{(out)} &  &  &  &  & \mbox{\Large{O}} \\ 
w_2^{(in)}w_1^{(out)} & 0 & w_2^{(in)}w_3^{(out)} &  &  &  &  \\ 
& w_3^{(in)}w_2^{(out)} & \cdot & \cdot &  &  &  \\ 
&  & \cdot & \cdot & \cdot &  &  \\ 
&  &  & \cdot & \cdot & \cdot &  \\ 
&  &  &  & \cdot & \cdot & w_{n-1}^{(in)}w_n^{(out)} \\ 
\mbox{\Large{O}} &  &  &  &  & w_n^{(in)}w_{n-1}^{(out)}  &0
\end{array}
\right]%\in\R^{n\times n}
\end{equation}}
with $0<w_j^{(in)},w_j^{(out)}\leq 1$ for all $j$; if person $k$ does not wear a mask, 
then $w_k^{(in)}=w_k^{(out)}=1$. This model assumes that all interactions are of the 
same duration and that the distance between adjacent people is the same; a rescaling of the 
$w_j^{(in)}$ and $w_j^{(out)}$ is required to model interactions of different durations
and of people being at different distances from each other. In any case, the matrix 
\eqref{Amat} typically is nonsymmetric.

We obtain an adjacency matrix that is simpler to analyze by projecting the matrix 
\eqref{Amat} orthogonally onto the subspace $\mathbb{T}$ of tridiagonal Toeplitz 
matrices of order $n$. Let $T$ be the orthogonal projection of the matrix \eqref{Amat} 
onto $\mathbb{T}$. Then
\begin{equation}\label{Tmath}
T=\left[ 
\begin{array}{ccccccc}
0 & t_1 &  &  &  &  & \mbox{\Large{O}} \\ 
t_{-1} & 0& t_1 &  &  &  &  \\ 
& t_{-1} & 0 & \cdot &  &  &  \\ 
&  & \cdot & \cdot & \cdot &  &  \\ 
&  &  & \cdot & \cdot & \cdot &  \\ 
&  &  &  & \cdot & \cdot & t_1 \\ 
\mbox{\Large{O}} &  &  &  &  & t_{-1} & 0
\end{array}
\right]\in\R^{n\times n},
\end{equation}
where the superdiagonal entry $t_1$ is the average of the superdiagonal entries of the
matrix \eqref{Amat}, and the subdiagonal entry $t_{-1}$ is the average of the subdiagonal 
entries of \eqref{Amat}; see, e.g., \cite{NPR07}. When all $w_j^{(in)}$ and $w_j^{(out)}$ are 
positive, so are $t_1$ and $t_{-1}$, and it follows that the matrix \eqref{Tmath} is
irreducible.

Assume for the moment that $w_j^{(in)}=w_j^{(out)}>0$ for all $j$. Then the matrices 
\eqref{Amat} and \eqref{Tmath} are symmetric. It follows from a result due to Bhatia 
\cite{RB86} that if the relative distance between the symmetric matrices \eqref{Amat} 
and \eqref{Tmath} is small in the Frobenius norm, then the relative difference of the 
spectra of \eqref{Amat} and \eqref{Tmath} also is small. In detail, let the matrices 
$M_1\in\R^{n\times n}$ and $M_2\in\R^{n\times n}$ be symmetric, and consider the 
relative distance between these matrices in the Frobenius norm, 
\[
d_{M_1,M_2}=\frac{\|M_1-M_2\|_F}{\|M_1\|_F}.
\]
Order the eigenvalues $\lambda_j(M_1)$ of $M_1$ and $\lambda_j(M_2)$ of $M_2$ according to 
$\lambda_1(M_1)\geq\lambda_2(M_1)\geq\ldots\geq\lambda_n(M_1)$ and
$\lambda_1(M_2)\geq\lambda_2(M_2)\geq\ldots\geq\lambda_n(M_2)$. Then
\[
\frac{\sqrt{\sum_{i=1}^{n}{(\lambda_i(M_1)-\lambda_i(M_2))^2}}}
{\sqrt{\sum_{i=1}^{n}{\lambda_i(M_1)^2}}}\leq d_{M_1,M_2}.
\]

However, as the following example shows, the spectral radius of $M_1$ may be much smaller 
than the spectral radius of $M_2$, also when $d_{M_1,M_2}$ is small. 

\begin{figure}
\centerline{
\includegraphics[scale=0.35,trim= 0mm 0.01mm 0mm 0mm]{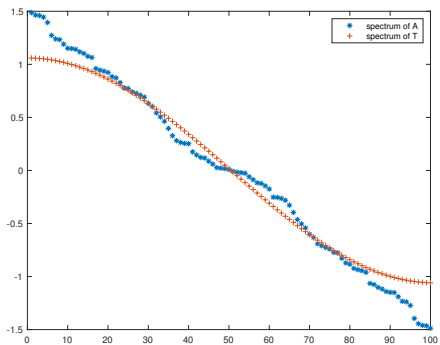}
}
\caption{Example \ref{ex1}: Eigenvalues of tridiagonal matrix $A$ (blue circle) and 
eigenvalues of the closest tridiagonal Toeplitz matrix $T$ (red circles).}\label{fig_ex}
\end{figure}

\begin{exmp}\label{ex1}
Let $A\in\R^{100\times 100}$ be a symmetric tridiagonal irreducible matrix with uniformly
distributed random entries in the interval $[0,1]$. These entries were generated with the
random number generator {\sf rand} in MATLAB. Let $T\in\R^{100\times 100}$ denote the 
closest symmetric tridiagonal Toeplitz matrix to $A$. We obtain $d_{A,T}=0.49$ and 
\[
\frac{\sqrt{\sum_{i=1}^{100}{(\lambda_i(A)-\lambda_i(T))^2}}}
{\sqrt{\sum_{i=1}^{100}{\lambda_i(A)^2}}}=0.19,
\]
where the eigenvalues of $A$ and $T$ are ordered in non-increasing order. Figure 
\ref{fig_ex} shows the eigenvalues of $A$ and $T$ as functions of their index. The extreme
eigenvalues of $A$ and $T$ are seen not to be close. In particular, the spectral radius of 
$T$ is quite a bit smaller than the spectral radius of $A$.
\end{exmp}

The following result is an analogue of Proposition \ref{peig} for nonsymmetric tridiagonal
Toeplitz matrices. 

\begin{proposition}\label{peig2}
The Perron root $\rho$ of the nonnegative tridiagonal Toeplitz matrix \eqref{Tmath} is 
$2\sqrt{t_{-1}t_1} \cos\frac{\pi}{n+1}$. The right Perron vector 
$\mathbf{u}=[u_1,u_2,\ldots,u_n]^T$, suitably scaled, has the entries 
$u_k=(t_{-1}/t_1)^{k/2}\sin \frac{k\pi}{n+1}$, $1\leq k\leq n$. The left Perron vector 
$\mathbf{v}=[v_1,v_2,\ldots,v_n]^T$, suitably scaled, has the entries 
$v_k=(t_1/t_{-1})^{k/2}\sin \frac{k\pi}{n+1}$, $1\leq k\leq n$.
\end{proposition}

\begin{proof}
Explicit formulas for eigenvalues and eigenvectors of tridiagonal Toeplitz matrices can be
found in, e.g., \cite{NPR}. 
\end{proof}

\begin{rmk}
Symmetrizing the matrix \eqref{Tmath}, i.e., considering an undirected graph instead of 
the directed graph represented by the adjacency matrix \eqref{Tmath}, gives the symmetric 
adjacency matrix $A=\frac{1}{2}(T+T^T)$ with Perron root 
\[
\rho(A)= 2\left(\frac{t_{-1}+t_1}{2}\right)\cos\frac{\pi}{n+1}\,. 
\]
Thus, the Perron root is determined by the arithmetic mean of $t_{-1}$ and $t_1$, while 
the Perron root of the matrix \eqref{Tmath} is defined by the geometric mean of these
quantities; cf. Proposition \ref{peig2}.
\end{rmk}

\begin{exmp}
Let $A \in \R^{25\times 25}$ be the  tridiagonal Toeplitz matrix \eqref{Tmath} with 
$t_{-1}=1.5$ and $t_{1}=0.5$.  This matrix may model a situation where the probability of 
inhaling infected droplets is three times larger than the probability of exhaling them; 
e.g., people wearing chirurgical masks.  Proposition \ref{peig2} yields $\rho(A)=1.719422$
and the unit norm right and left Perron vectors $\mathbf{u}$ and $\mathbf{v}$. It is easy 
to see that edge $e_{12,13}$ is a maximizer of $\max_{h,k} \alpha_{hk}$, where 
$\alpha_{hk}=a_{hk}v_hu_k$; see \eqref{sp_rad}--\eqref{alphacoeff}.
In order to reduce the weight of $e_{12,13}$, one constructs the perturbed matrix
\[
\widetilde{A}=A+\varepsilon E_{12,13}= A-\varepsilon a_{12,13}\mathbf{e}_{12}
\mathbf{e}_{13}^T.
\]
Setting $\varepsilon=0.1$, we obtain $\rho(\widetilde{A})=1.713348$. Thus, the spectral impact 
of the perturbation is $s_{12,13}^{(\rho(A))}=0.003532$; its first order approximation 
\[
\alpha_{12,13}\frac{\varepsilon\kappa(\rho(A))}{\rho(A)}=\,a_{12,13}v_{12} u_{13}
\frac{\varepsilon\kappa(\rho(A))}{\rho(A)}=0.003846
\]
is fairly close. 

Assume that there is only one high-quality face mask available. This example shows which person 
should be wearing it to reduce the spectral radius the most. Notice that symmetrizing the 
matrix $A$ would have given both the matrix and the results of Example \ref{tridtoep}.
\end{exmp}

\section{Estimating and reducing the spectral radius}\label{sec3}
This section discusses several ways to estimate the spectral radius, and the right and 
left Perron vectors, of a large adjacency matrix $A\in\R^{n\times n}$. If $A$ just is
required to be nonnegative, then there is a nonnegative vector $\mathbf{x}\in\R^n$, such 
that $A\mathbf{x}=\rho(A)\mathbf{x}$. However, this vector may not be unique up 
to scaling; see \cite[Theorem 8.3.1 and p. 505]{HJ}. In this section, we will assume that 
$A$ is a nonnegative irreducible adjacency matrix. Then its right and left Perron vectors 
are unique up to scaling, and can be scaled to be of unit norm and have positive entries
only. These vectors
are used to determine which edge-weights to reduce to obtain a new adjacency matrix with,
hopefully, a significantly reduced spectral radius. If our aim just is to determine the
spectral radius of $A$, then irreducibility is not required.

Moreover, we analyze the behavior of the Perron root  when the adjacency matrix is 
perturbed by the worst $\varepsilon$-size perturbation for $\rho(A)$. This approach of 
approximating the rightmost $\varepsilon$-pseudoeigenvalue of $A$ takes into account 
the possibility of the entries of $A$ being contaminated by error and gives an 
estimate of the sensitivity of the structural robustness to a worst-case perturbation.

We first describe a computational method that is well suited for large networks, whose 
associated adjacency matrix is nonnegative and irreducible, but does not have other 
structure that can be exploited. Subsequently, we will discuss methods that are able to 
use certain structural properties.

\subsection{Approximation of the spectral radius of a nonnegative irreducible 
matrix}\label{alg1}
Let $A\in\R^{n\times n}$ be a large nonnegative irreducible adjacency matrix. The approach
of this section does not exploit any additional structure that $A$ may possess. We 
determine approximations of the right and left Perron vectors of $A$ by the two-sided 
Arnoldi method. This method was first described by Ruhe \cite{Ru} and has more recently 
been studied and improved by Zwaan and Hochstenbach \cite{ZH}. 

We carry out the following steps:
\begin{itemize}
\item
Apply the two-sided Arnoldi method to $A$ to compute the Perron root $\rho(A)$, and the 
unit right and left Perron vectors $\mathbf{u}$ and $\mathbf{v}$, respectively, with
positive entries.
\item
Let 
\begin{equation}\label{matE}
E=\mathbf{v}\mathbf{u}^T.
\end{equation}
For each edge $e_{hk}\in\cE$ in the graph that represents the network, the corresponding 
entry of $E$, i.e., $v_hu_k$, appears in the first-order approximation \eqref{alphacoeff} 
of the spectral impact of the edge. The Perron root $\rho(A+\varepsilon E)$ of the matrix 
$A+\varepsilon E$ satisfies
\[
\rho(A+\varepsilon E)=\rho(A)+\varepsilon 
\frac{\mathbf{v}^TE\mathbf{u}}{\mathbf{v}^T\mathbf{u}}+O(\varepsilon^2)
\]
for $|\varepsilon|$ sufficiently small; see Wilkinson \cite[Chapter 2]{Wi}. We refer to the
matrix \eqref{matE} as a \emph{Wilkinson perturbation}. This is the worst perturbation 
for $\rho(A)$ in the following sense. For any nonnegative matrix $E$ with $\|E\|=1$, one has
\[
\frac{\mathbf{v}^T E \mathbf{u}} {\mathbf{v}^T\mathbf{u}}=
\frac{|\mathbf{v}^T E \mathbf{u}|} {\mathbf{v}^T\mathbf{u}}\leq
\frac{\|\mathbf{v}\|\|E\|\|\mathbf{u}\|} {\mathbf{v}^T\mathbf{u}}=
\frac{1}{\mathbf{v}^T\mathbf{u}},
\]
with equality for the matrix \eqref{matE}. Moreover,
$$
\rho(A+\varepsilon E)-\rho(A)\approx \varepsilon\frac{\mathbf{v}^T E \mathbf{u}} {\mathbf{v}^T\mathbf{u}}
=\varepsilon \kappa(\rho(A)).
$$
We let $\varepsilon>0$. 
Note that the spectrum of $A+\varepsilon E$ may be considered a very sparse approximation 
of the $\varepsilon$-pseudospectrum of $A$ in the sense that the 
$\varepsilon$-pseudospectrum is made up of the spectra of all perturbations of $A$ of
norm at most $\varepsilon$; see Trefethen and Embree \cite{TE}. We only compute the
spectrum for one perturbation, $A+\varepsilon E$, of $A$. The size of $\varepsilon$ used 
in the computations may depend on whether the adjacency matrix is contaminated by errors. 
For instance, the edge weights may not be known exactly; see Trefethen and Embree 
\cite{TE} for insightful discussions on pseudospectra.
\item
Typically, the first order approximation 
\[
\rho(A)+\varepsilon\frac{\mathbf{v}^T E \mathbf{u}} {\mathbf{v}^T\mathbf{u}}
\]
of $\rho(A+\varepsilon E)$ is sufficiently accurate. In the rare occasions when this is not
the case, we can compute an improved approximation by applying the (standard) Arnoldi
method described, e.g., by Saad \cite{Sa}, or the implicitly restarted (standard) Arnoldi
method described in \cite{LSY} and implemented by the MATLAB function {\sf eigs}.
We note that the perturbed matrix $A+\varepsilon E$ is nonnegative and irreducible if this
holds for $A$. Indeed, if all entries of the Perron vectors are positive, then so are all 
entries of $A+\varepsilon E$ for $\varepsilon>0$.
\end{itemize}

The Perron root $\rho(A+\varepsilon E)$ is a rightmost $\varepsilon$-pseudoeigenvalue of 
$A$. We note that $\rho(A+\varepsilon E)$ may be much larger than $\rho(A)$ when the 
Perron root is ill-conditioned, i.e., when $\mathbf{v}^T\mathbf{u}$ is small.

The analysis in Section \ref{sec2} suggests that in order to reduce $\rho(A)$ by removing
an edge of $\cG$, we should choose an edge $e_{hk}$ with a large weight $a_{hk}$ that 
corresponds to a large entry of the matrix $\mathbf{v}\mathbf{u}^T$ in \eqref{matE}; see 
\eqref{sp_rad}-\eqref{alphacoeff}. Removing an edge corresponds to setting its 
edge-weight to zero. We can in the same manner choose which edge-weight to reduce to a 
smaller positive value in order to reduce the spectral radius.

\begin{figure}
\centerline{
\includegraphics[scale=0.35,trim= 0mm 0.01mm 0mm 0mm]{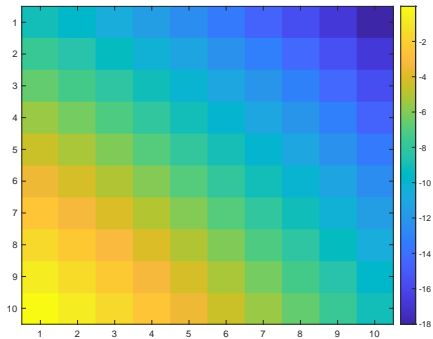}
}
\caption{Example \ref{circ}: $\log_{10}$ of the entries of the matrix \eqref{matE}.} 
\label{fig_circ}
\end{figure}

\begin{exmp} \label{circ}
Consider a matrix $T\in\mathbb{R}^{10\times 10}$ of the form \eqref{Tmath}
with $t_{-1}=0.1$ and $t_1=1$. The eigenvalues  of $T$ are real and appear in $\pm$ pairs.
Thus, there are two eigenvalues of largest magnitude. The positive one is about $0.6$. Now
we add suitable entries in the $(1,10)$ and $(10,1)$ positions to transform $T$ into a
circulant matrix $T_c$. Then $\rho(T_c)=1.1$ and $T_c$ also has the eigenvalue $-1.1$. The 
remaining eigenvalues are complex-valued.  

The large perturbation induced in the spectrum of $T$ by this perturbation of $T$ can be 
explained  by analyzing the structure of the matrix $E$ in \eqref{matE}, which we 
construct by using the left and right Perron vectors given in Proposition \ref{peig2}. 
Figure \ref{fig_circ} visualizes the size of the entries in $E$. Notice that the largest 
entries are confined to the bottom left corner. Thus, adding the entry $1$ in the 
position $(10,1)$ induces a large perturbation in the Perron root.
\end{exmp}

When the adjacency matrix $A$ is very large, we may consider replacing the vectors 
$\mathbf{u}$ and $\mathbf{v}$ in \eqref{matE} by the vector 
$\mathbf{e}=\frac{1}{\sqrt{n}}[1,1,\ldots,1]^T$ and compute $\rho(A)$ and 
$\rho(A+\varepsilon\, \mathbf{e}\,\mathbf{e}^T)$ by the (standard) Arnoldi or restarted 
Arnoldi methods to determine the structural robustness of the graph with adjacency matrix
$A$. This approach was applied in \cite{NR} to estimate pseudospectra of large matrices.

The large perturbation in $\rho(A)$ illustrated in Example \ref{circ} would not have 
occurred if the sparsity structure of the matrix $T$ would have been taken into account,
i.e., if one only would allow perturbations of positive edge-weights. We therefore are
interested in determining perturbations $\varepsilon E$ of $A$ that take the sparsity 
structure of $A$ into account.

\subsection{Approximation of the spectral radius taking the sparsity structure into 
account}\label{alg2}
The method in this subsection is suitable when it is desirable that the perturbation 
$\varepsilon E$ of the adjacency matrix $A$ has the same sparsity structure as $A$. Let
${\mathcal S}$ denote the cone of all nonnegative matrices in ${\mathbb R}^{n\times n}$
with same sparsity structure as $A$, and let $M|_{\mathcal S}$ be the matrix in ${\mathcal S}$ 
that is closest to a given nonnegative matrix $M$ with respect to the Frobenius norm. 
It is straightforward to verify that
the matrix $M|_{\mathcal S}$ is obtained by replacing all the entries of $M$ outside the 
sparsity structure by zero. This approach takes possible uncertainty of the available 
edge-weights into account. The analysis in \cite{NR} leads to the following numerical 
method:
\begin{itemize}
\item
Apply the two-sided Arnoldi method to $A\in{\mathcal S}$ to compute the Perron root 
$\rho(A)$, as well as the unit right and left Perron vectors $\mathbf{u}$ and 
$\mathbf{v}$, respectively, with positive entries.
\item 
Project $\mathbf{v}\mathbf{u}^T$ into ${\mathcal S}$,
normalize the projected matrix to have unit Frobenius norm, and define
\begin{equation}\label{matES}
E=\frac{\mathbf{v}\mathbf{u}^T|_{\mathcal S}}{\| \mathbf{v}\mathbf{u}^T|_{\mathcal S} \|_F}.
\end{equation}
We refer to the matrix \eqref{matES} as an \emph{${\mathcal S}$}-structured analogue of the 
Wilkinson perturbation. This is the worst ${\mathcal S}$-structured perturbation 
for $\rho(A)$; one has, by \cite[Proposition 2.3]{NP06},
\[
\frac{\mathbf{v}^TE\mathbf{u}}{\mathbf{v}^T\mathbf{u}}=
\frac{|\mathbf{v}^TE\mathbf{u}|}{\mathbf{v}^T\mathbf{u}}=
\frac{\|\mathbf{v}\|\|\mathbf{v}\mathbf{u}^T|_{\mathcal S}\|_F\|\mathbf{u}\|}
{\mathbf{v}^T\mathbf{u}}=
\frac{\|\mathbf{v}\mathbf{u}^T|_{\mathcal S}\|_F}{\mathbf{v}^T\mathbf{u}}.
\]
Hence,
$$
\rho(A+\varepsilon E)-\rho(A)\approx 
\varepsilon\frac{\mathbf{v}^T E \mathbf{u}} {\mathbf{v}^T\mathbf{u}}
= \varepsilon \kappa_{\mathcal S}(\rho(A)), 
$$
where
\[
\kappa_{\mathcal S}(\rho(A))= \frac{\| \mathbf{v}\mathbf{u}^T|_{\mathcal S} \|_F}
{\mathbf{v} ^T \mathbf{u}}
\]
denotes the \emph{${\mathcal S}$-structured condition number} of $\rho(A)$; see 
\cite{NP06,KKT}. We let $\varepsilon>0$. Similarly as above, the spectrum of $A+\varepsilon E$
is a very sparse approximation of the ${\mathcal S}$-structured 
$\varepsilon$-pseudospectrum of $A$ in the sense that this pseudospectrum is evaluated by
computing the spectrum for many ${\mathcal S}$-structured perturbations of $T$ of norm at 
most $\varepsilon$; see, e.g., \cite{NR}. Here we only compute the spectrum for one 
perturbation of $T$.
\item
If desired, compute $\rho(A+\varepsilon E)$ by the (standard) Arnoldi or restarted Arnoldi
methods. We note that the perturbed matrix $A+\varepsilon E$ is nonnegative and 
irreducible if this holds for $A$, and exhibits the same sparsity structure as $A$.
\end{itemize}

The Perron root $\rho(A+\varepsilon E)$ helps us to estimate the structural robustness of
the network. Indeed, it represents an approximate ${\mathcal S}$-structured 
$\varepsilon$-pseudospectral radius of the ${\mathcal S}$-structured 
$\varepsilon$-pseudospectrum of the adjacency matrix $A \in {\mathcal S}$.

We note that $\rho(A+\varepsilon E)$ may be much 
larger than $\rho(A)$ when the Perron root has a large  ${\mathcal S}$-structured 
condition number $\kappa_{\mathcal S}(\rho(A))$. 

As mentioned above, in case the network is very large, we may consider replacing the 
vectors $\mathbf{u}$ and $\mathbf{v}$ in \eqref{matES} by the vector $\mathbf{e}$. An
analogous \emph{${\mathcal S}$-structured perturbation} of the adjacency matrix $A$ is 
given by 
\[
\frac{\mathbf{e}\mathbf{e}^T|_{\mathcal S}}{\|\mathbf{e}\mathbf{e}^T|_{\mathcal S}\|_F}.
\]
We may apply the (standard) Arnoldi or 
implicitly restarted Arnoldi methods to estimate $\rho(A)$ and 
\[
\rho\left(A+\varepsilon\, \frac{\mathbf{e}\mathbf{e}^T|_{\mathcal S}}
{\| \mathbf{e}\mathbf{e}^T|_{\mathcal S} \|_F}\right).
\]
This approach has been applied in \cite{NR} to estimate structured pseudospectra of large 
matrices.

\subsection{Approximation of the spectral radius for perturbations of tridiagonal Toeplitz
matrices}\label{alg3}
Structure respecting projections, analogous to the ones discussed in the previous subsection, 
also can be applied to impose other structures. This subsection illustrates how they can 
be used to impose tridiagonal Toeplitz structure. Let $T$ be a nonnegative tridiagonal
Toeplitz matrix \eqref{Tmath}. We denote by ${\mathcal T}$ the cone of all nonnegative
tridiagonal Toeplitz matrices with zero diagonal in  ${\mathbb R}^{n\times n}$ and by 
$M|_{{\mathcal T}}$ the matrix in ${\mathcal T}$ closest to a given nonnegative matrix
$M\in {\mathbb R}^{n\times n}$
with respect to the Frobenius norm. It is straightforward to verify that 
$M|_{{\mathcal T}}$ is obtained by replacing the sub- and super- diagonal entries of $M$ 
by their respective arithmetic mean. 

To approximate the spectral radius of $T\in{\mathcal T}$, we carry out the following 
steps:
\begin{itemize}
\item
Apply the formulas in Proposition \ref{peig2} to $T$ to compute the Perron root $\rho(T)$
and the unit right and left Perron vectors $\mathbf{u}$ and $\mathbf{v}$, respectively, 
with positive entries.
\item
Project $\mathbf{v}\mathbf{u}^T$ into ${\mathcal T}$, normalize the projected matrix
to have unit Frobenius norm, and define the matrix
\begin{equation}\label{matET}
E=\frac{\mathbf{v}\mathbf{u}^T|_{\mathcal T}}{\|\mathbf{v}\mathbf{u}^T|_{\mathcal T} \|_F}.
\end{equation}
We refer to the matrix \eqref{matET} as a \emph{${\mathcal T}$-structured analogue of the 
Wilkinson perturbation}. Similarly as above, we have, by \cite[Theorem 3.3]{NP07},
\[
\frac{\mathbf{v}^TE\mathbf{u}}{\mathbf{v}^T\mathbf{u}}=
\frac{|\mathbf{v}^TE\mathbf{u}|}{\mathbf{v}^T\mathbf{u}}=
\frac{\|\mathbf{v}\|\|\mathbf{v}\mathbf{u}^T|_{\mathcal T}\|_F\|\mathbf{u}\|}{\mathbf{v}^T\mathbf{u}}=
\frac{\|\mathbf{v}\mathbf{u}^T|_{\mathcal T}\|_F}{\mathbf{v}^T\mathbf{u}}.
\]
It follows that
$$
\rho(T+\varepsilon E)-\rho(T)\approx\varepsilon 
\frac{\mathbf{v}^TE\mathbf{u}}{\mathbf{v}^T\mathbf{u}}=
\varepsilon \kappa_{\mathcal T}(\rho(T)),
$$
where
\[
\kappa_{\mathcal T}(\rho(T))=
\frac{\|\mathbf{v}\mathbf{u}^T|_{\mathcal T} \|_F}{\mathbf{v}^T \mathbf{u}}
\]
denotes the \emph {${\mathcal T}$-structured condition number} of $\rho(T)$; see 
\cite{NP07,KKT}. 

We will let $\varepsilon>0$. The spectrum of $T+\varepsilon E$ is a sparse approximation 
of the ${\mathcal T}$-structured $\varepsilon$-pseudospectrum of $T$ in the sense that 
this pseudospectrum is evaluated by computing the spectrum for many 
${\mathcal T}$-structured perturbations of $T$ of norm at most $\varepsilon$; see, e.g.,
\cite{NR}. Here we only compute the spectrum for one perturbation of $T$.
\item
Determine $\rho(T+\varepsilon E)$ by applying Proposition \ref{peig2} to 
$T+\varepsilon E$. The latter matrix is nonnegative and irreducible if this holds for $T$,
and exhibits the same structure as $T$.
\end{itemize}

The Perron root $\rho(T+\varepsilon E)$ may be regarded as an approximate 
${\mathcal T}$-structured $\varepsilon$-pseudospectral radius and provides an estimate of
the structural robustness of the structured network. It may be much larger than $\rho(T)$.
It is known that when considering the class ${\mathcal T}$ of tridiagonal Toeplitz 
matrices, the most ill-conditioned eigenvalues with regard to ${\mathcal T}$-structured 
perturbations are the eigenvalues of largest magnitude; see, e.g., \cite{NPR}. We remark 
that an algorithm for computing the ${\mathcal T}$-structured pseudospectrum of a 
tridiagonal Toeplitz matrix and its rightmost pseudoeigenvalue is described in \cite{BGN}.
However, the computational cost of this algorithm can be quite large for the matrices
considered in this paper.

Finally, replacing $\mathbf{u}$ and $\mathbf{v}$ in \eqref{matET} by the vector 
$\mathbf{e}$ as described above is particularly efficient when the considered subspace is
${\mathcal T}$; see Section \ref{test_tridtoep}.

\section{Numerical tests}\label{sec4}
his section illustrates the performance of the methods discussed when applied to large
networks. All computations were carried out in MATLAB R2021a on a MacBook Pro with a
2GHz Intel Core i5 quad-core CPU and 16GB of RAM.

\subsection{Complex networks}\label{comnet}

\subsubsection{Air500}\label{test1}
Consider the adjacency matrix $A\in \R^{500 \times 500}$ for the network Air500, which
describes 24009 flight connections between the top $500$ airports within the United States
based on total passenger volume during one year from July 1, 2007, to June 30, 2008; 
see~\cite{Bio}. Thus, the airports are nodes and the flights are edges in the graph 
determined by the network. The matrix $A$ has the entry $a_{i,j}=1$ if there is a flight 
that leaves from airport $i$ to airport~$j$. Generally, but not always, $a_{i,j}=1$ 
implies that $a_{j,i}=1$. This makes $A$ close to symmetric.

\begin{figure}
\centerline{
\includegraphics[scale=0.35,trim= 0mm 0.01mm 0mm 0mm]{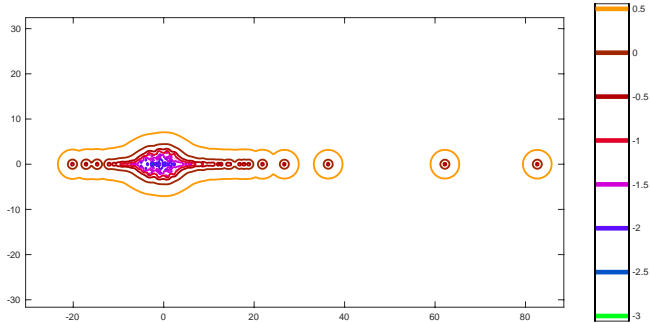}
}
\caption{Pseudospectra of the adjacency matrix of the network in Section \ref{test1}. The
contour levels are from $10^{-3}$ to $10^{0.5}$. The figure is produced by EigTool.} 
\label{fig_test1}
\end{figure}

Apply the computational steps described in Section \ref{alg1}. The Perron root $\rho(A)$ 
is $82.610276$ with eigenvalue condition number $\kappa(\rho(A))=1.001668$. Let
$\varepsilon = 0.5$. The Perron root $\rho(A+\varepsilon E)$, where $E$ is the matrix in
\eqref{matE}, is $83.111096$. Thus, the spectral radius increases by $0.500820$, as we 
could have foreseen since $\varepsilon\kappa(\rho(A))=0.500834$. The value 
$\rho(A+\varepsilon E)$ is an accurate approximation of the $\varepsilon$-pseudospectral 
radius. This is seen by determining the $\varepsilon$-pseudospectral radius by the MATLAB
program package Eigtool \cite{Wr}. Our approximation of the $\varepsilon$-pseudospectral 
radius agrees with the value determined by Eigtool in all decimal digits returned by 
Eigtool. Pseudospectra of $A$ are visualized in Figure \ref{fig_test1}.

Assume that we are interested in removing a single route so that the structural robustness 
of the network is increased the most. Then this route should be an edge that maximizes 
$a_{h,k} v_hu_k$ over $h$ and $k$; see \eqref{sp_rad}-\eqref{alphacoeff}. For the present
network, we find that the edge $e_{224,24}\in \mathcal{E}$ should be removed. The adjacency
matrix $\widehat{A}$ so obtained is irreducible with $\rho(\widehat{A})=82.590199$. The 
edge $e_{224,24}$ corresponds to flights from the JFK airport in New York to the 
Hartsfield--Jackson airport in Atlanta.

Finally, we observe that if one replaces $E$ in \eqref{matE} by the
matrix of all ones normalized to have unit Frobenius norm, the  increase of spectral 
radius $\rho(A)$ results to be $0.255450$. Thus, this perturbation gives a significantly
less accurate estimate of the sensitivity of $\rho(A)$ to worst-case perturbations.

\subsubsection{Airlines}\label{test2}
Consider the adjacency matrix $A\in \R^{235 \times 235}$ determined by the network 
Airlines with 235 nodes and 2101 edges. The nodes represent airports and the directed 
edges represent flights between them. This network is available at \cite{Ge}. 

\begin{figure}
\centerline{
\includegraphics[scale=0.35,trim= 0mm 0.01mm 0mm 0mm]{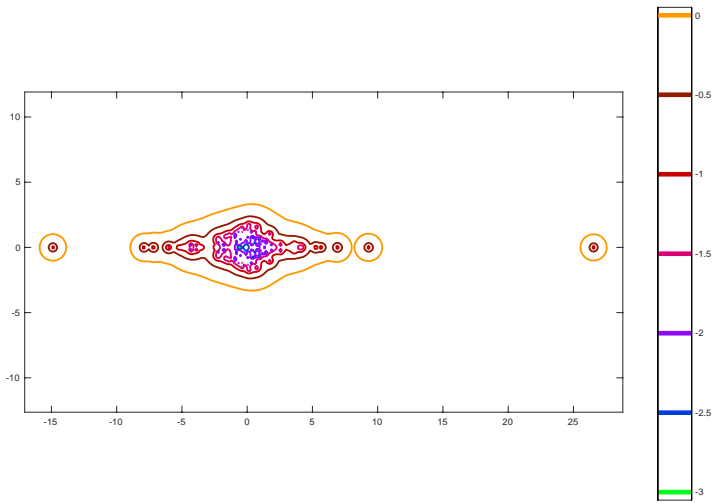}
}
\caption{Pseudospectra of the adjacency matrix of the network in Section \ref{test2}. 
Contour curves are for levels from $10^{-3}$ to $10^{0}$. The graph is produced by EigTool.}
\label{fig_test2}
\end{figure}

Computations described in Section \ref{alg1} yield $\rho(A)=26.545430$ and the 
condition number $\kappa(\rho(A))=1.005219$. Let $\varepsilon = 0.5$. The Perron root 
$\rho(A+\varepsilon E)$, where $E$ is the matrix \eqref{matE}, is $27.047941$. Thus, 
the spectral radius increases by $0.502511$, as we could have expected since 
$\varepsilon\kappa(\rho(A))=0.502609$. The spectral radius 
$\rho(A+\varepsilon E)$ approximates the $\varepsilon$-pseudospectral radius 
and agrees with six significant decimal digits with the pseudospectral radius 
determined by Eigtool. A few pseudospectra of $A$ are shown in Figure \ref{fig_test2}.

The route to remove, in order to increase the structural robustness of the network the 
most, is represented by the edge $e_{51,137}\in \mathcal{E}$. The adjacency matrix, 
$\widehat{A}$, obtained when setting the entry $a_{51,137}$ of $A$ to zero is irreducible 
with $\rho(\widehat{A})=26.452922$.

Finally, we observe that if one replaces $E$ in \eqref{matE} by the matrix of all ones,
normalized to have unit Frobenius norm, $\rho(A)$ increases by $0.223135$. 

\subsubsection{Enron}\label{test3} 
The Enron e-mail exchange network covers the e-mail communication at the Enron company.
The data set consists of 
over $3\cdot 10^5$ e-mails. The e-mail addresses are the nodes $n_i$ of the network; 
there are 36692 of them. A directed edge from node $n_i$ to node $n_j$ indicates that at 
least one e-mail message was sent from $n_i$ to $n_j$; there are 367662 edges. Let 
$A\in \R^{36692 \times 36692}$ be the adjacency matrix for this graph. It is close to 
symmetric. This network is available at \cite{SN}. 

\begin{figure}
\centerline{
\includegraphics[scale=0.5,trim= 0mm 0.01mm 0mm 0mm]{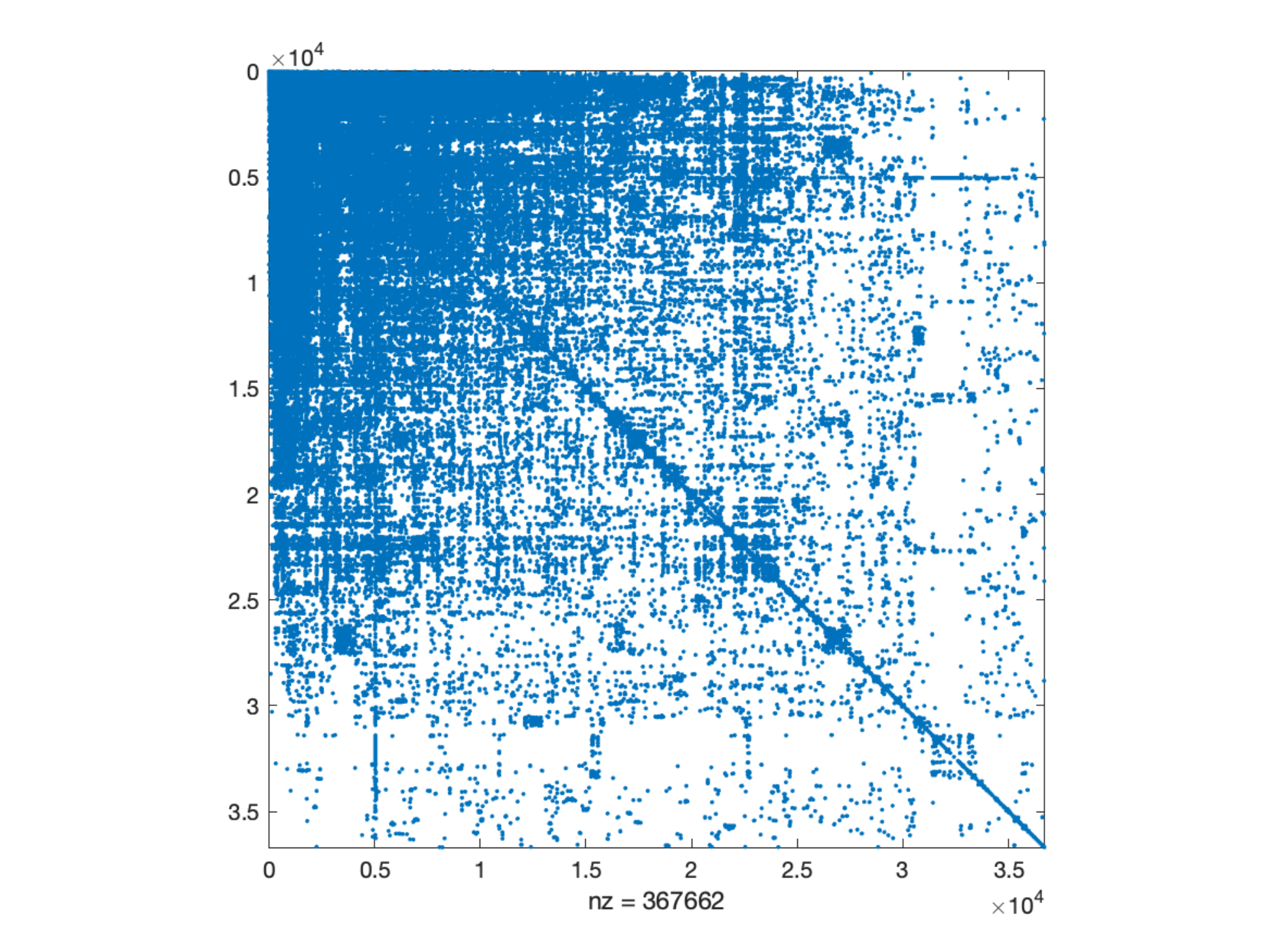}
}
\caption{Spy plot of the adjacency matrix of the network in Section \ref{test3}.}
\label{fig_test3}
\end{figure}

Computations described in Section \ref{alg1} yield $\rho(A)=118.417715$ and the condition 
number $\kappa(\rho(A))=1.000000$. Let $\varepsilon = 0.5$. 
The Perron root is $\rho(A+\varepsilon E)=118.917715$, where $E$ is the matrix 
\eqref{matE}. As expected, the spectral radius increases by $0.500000$. No comparison 
with Eigtool could be made to assess the accuracy of the so determined approximation
of the $\varepsilon$-pseudospectral radius, since Eigtool is not able to determine
$\varepsilon$-pseudospectra of such a large matrix. The spy plot of $A$ is shown in 
Figure \ref{fig_test3}. In this plot positive matrix entries are marked by a dot.

The email-channel to remove is represented by the edge $e_{137,196}\in \mathcal{E}$. As 
for the adjacency matrix $\widehat{A}$, obtained when setting the entry $a_{137,196}$ of
$A$ to zero, one has $\rho(\widehat{A})=118.398705$. In fact, it is immediate to see from 
the spy plot in Figure \ref{fig_test3} that the edge to be removed belongs to a ``dense 
block'' of $A$. Removing this edge yields a relative decrease of order $10^{-4}$ in 
the spectral radius of the adjacency matrix of the Enron network.

\subsection{Synthetic networks}
This subsection considers projections of the adjacency matrix for the Air500 network.

\subsubsection{The tridiagonal part of Air500}\label{test_trid}
We set all entries of the adjacency matrix for the Air500 network outside the tridiagonal 
part of the adjacency matrix to zero. The number of flight connections is now $144$.
This yields a nonsymmetric tridiagonal matrix
$A\in \R^{500 \times 500}$. Carry out the computations described in Section \ref{alg2}, 
with ${\mathcal S}$ the subspace of all tridiagonal matrices with zero-diagonal in 
${\mathbb R}^{500\times 500}$. This yields the Perron root $\rho(A)=1.801938$ and its 
${\mathcal S}$-structured condition number is $\kappa_{\mathcal S}(\rho(A))=0.613714$.

Let $\varepsilon = 0.9$. The Perron root $\rho(A+\varepsilon E)$, where $E$ is 
the matrix in \eqref{matES}, is $2.362116$. Thus, the spectral radius increases by 
$0.560178$, as we could have foreseen since 
$\varepsilon\kappa_{\mathcal S}(\rho(A))=0.552342$. 

Computations similar to those of Subsection \ref{comnet} suggest that in order to increase
the structural stability the most by removing one edge, we should choose the edge 
$e_{494,493}$ or the edge $e_{493,494}$ in $\mathcal{E}$. However, removal of one or both of these 
edges would result in a graph with a reducible adjacency matrix. To preserve 
irreducibility of the adjacency matrix, one may instead schedule fewer flights on the 
routes that correspond to the edges $e_{494,493}$ and $e_{493,494}$. This reduces the 
weight associated with these edges.

Finally, we observe that, if one replaces the matrix $\mathbf{v}\mathbf{u}^T$ in 
\eqref{matES} by the matrix of all ones, normalized to be of unit Frobenius norm, then the
spectral radius increases by $0.052606$. Clearly, this is not an accurate estimate of the
actual worst-case sensitivity of $\rho(A)$ to perturbations.

\subsubsection{Projection of Air500 into a tridiagonal Toeplitz structure}
\label{test_tridtoep}
We construct a tridiagonal Toeplitz matrix with zero-diagonal $T\in\R^{500 \times 500}$
by averaging the sub-diagonal entries as well as by averaging the super-diagonals of the 
matrix in Section \ref{test_trid}. Then 
we carry out the computations as described in Section \ref{alg3}, and make use of 
Proposition \ref{peig2}. We obtain $\rho(T)=0.288460$ and the ${\mathcal T}$-structured 
condition number $\kappa_{\mathcal T}(\rho(T))=0.063357$.

Let $\varepsilon=0.9$. Then $\rho(T+\varepsilon E)=0.345466$, where $E$ is 
the matrix in \eqref{matET}. Thus, the spectral radius increases by $0.057006$. This is in
agreement with $\varepsilon\kappa_{\mathcal T}(\rho(T))=0.057021$. 

Finally, we observe that, if one replaces the matrix $\mathbf{v}\mathbf{u}^T$ in 
\eqref{matES} by the matrix of all ones, scaled to be of unit Frobenius norm, then 
$\rho(T)$ increases by $0.056995$. Thus, the latter perturbation provides a very 
accurate estimate of the spectral radius when the matrix $\mathbf{v}\mathbf{u}^T$ in
\eqref{matET} is used. 

\section{Conclusion}\label{sec5}
It is important to be able to estimate the structural robustness of a network, and to
determine which nodes to remove or weights to decrease to increase the structural 
robustness. This paper describes several iterative methods that can be applied to fairly
large networks to gain insight into these issues. Both the sensitivity of the 
structural robustness to worst-case Wilkinson perturbations and to structured
perturbations are discussed and illustrated.

\section*{Acknowledgment}
The authors would like to thank Ian Zwaan for MATLAB code for the two-sided Arnoldi method
used in the numerical experiments and the anonymous reviewers for providing insightful 
comments. Work by SN was partially supported by a grant from SAPIENZA Universit\`a di Roma and by INdAM-GNCS, and work by LR was supported in part by NSF grant DMS-1720259. 

%\section*{Data availability statement}
%Data sharing not applicable to this article as no datasets were generated during the current study.

\end{document}